\documentclass{svmult-turin}

\usepackage{type1cm} 

\usepackage{makeidx}         
\usepackage{graphicx}        

\usepackage[bottom]{footmisc}

\usepackage{newtxtext}      %
\usepackage{newtxmath}      

\numberwithin{equation}{section}

\def\proof{\noindent\hspace{2em}{\itshape Proof: }}
\def\QEDclosed{\mbox{\rule[0pt]{1.3ex}{1.3ex}}} 

\def\QED{\QEDclosed} 
\def\endproof{\hspace*{\fill}~\QED\par\endtrivlist\unskip}
\newcommand{\eqa}{\begin{eqnarray}}
\newcommand{\eeqa}{\end{eqnarray}}
\newcommand{\beq}{\begin{equation}}
\newcommand{\eeq}{\end{equation}}

\begin{document}

\title*{On near-optimal time samplings for initial data best approximation}
\titlerunning{On near-optimal time samplings for initial data best approximation}

\author{Roza Aceska, Alessandro Arsie  and Ramesh Karki}
 
%
%
\maketitle

\abstract{ 
Leveraging on the work of De Vore and Zuazua, we further explore their methodology and deal with two open questions presented in their paper. We show that for a class of linear evolutionary PDEs  the admissible choice of relevant parameters used to construct the near-optimal sampling sequence is not influenced by the spectrum of of the operator controlling the spatial part of the PDE, but only by its order. Furthermore, we show that it is possible to extend their algorithm to a simple version of a non-autonomous heat equation in which the heat diffusivity coefficient depends explicitly on time.}

\section{Introduction}
\label{JSsec:introduction}

 The determination or the best possible approximation of the initial state of a dynamical system through observation of its states at subsequent times is a general and very important problem for a variety of applications. In \cite{DZ}, the authors proposed an ingenious procedure to approximate in a near-optimal way and via finitely many time samplings at a fixed location the initial state of a particularly simple infinite dimensional dynamical system, described by the heat equation on a compact interval with Dirichlet boundary conditions. 
 
We study an initial value problem for two generalized classes of PDEs,  involving fairly unknown initial conditions; each of these classes contain the heat equation as an example. It is known that under appropriate assumptions  (\cite{DZ}), one can compensate for the lack of knowledge of the initial condition by adding scarce measurements made at later time instances. This problem of compensation via a time-space trade off between the initial measurements and the later time measurements has been recently observed in applications of sampling theory, and referred to as the dynamical sampling problem (see for instance \cite{ACMT}, \cite{AP}, \cite{ADK12}).

We generalize the method developed in  \cite{DZ} and deal with an inverse problem (initial data best reconstruction) via later time measurements. The contributions of our work are significant in applications, where full knowledge of the initial conditions is unrealistic to expect. We study the correlation between the number of measurements that are needed to recover the initial profile to a prescribed accuracy, give precise estimates for the time instances when these measurements need to occur, and provide an optimal reconstruction algorithm under the assumption that the initial profile is in a Sobolev class. Let us underline that the inverse problem for the PDEs we are dealing with is in general ill-posed (see for instance \cite{BBC}, \cite{C}), however, our goal here is simply to find a best approximation of the initial data in $L^2([0,\pi])$ using finitely many time samplings at a fixed location. 

In  \cite{DZ}, the authors pointed out some open questions, some of which are dealt with in this paper.  
 
 The first question concerns the relationship between the spectrum of a certain operator and the choice of the parameter $\rho$ that determines the geometric sequence of near-optimal sampling times. In particular, the authors of \cite{DZ} inquired about the dependence of $\rho$ on the spectrum of a certain operator. In Section \ref{Section 2}, we study a generalization of the heat equation, essentially a constant coefficients evolutionary PDE of spatial order $2N$ with Dirichlet boundary conditions. In order the make our results more easily comparable with \cite{DZ}, we will make the strong assumption that the initial datum $f$ lives in a suitable subspace $S$ of the Sobolev space $H^{r}_0([0,\pi])$ for $r>0$. For our model, we will see that any $\rho>2N\ln(2)$ ($2N$ is the spatial order of the PDE)  will generate a near-optimal geometric sampling times sequence. In particular, this does not depend on the cofficients of the equation (as long as they satisfy a suitable sign requirement for well-posedness). So in this sense the dependence of $\rho$ on the spectrum is very weak. (In our case, the relevant spectrum is given by $\{\lambda(k)\}_{k\in \mathbb{N}},$ where $\lambda(k)=\sum_{l=1}^N (-1)^l \alpha_{2l}k^{2l}$, where $\alpha_{2l}$ are the constants appearing in the spatial part of the PDE.) Let us remark that in this Section the extension of the results of \cite{DZ}, although elementary, is not completely straightforward. 
In \cite{DZ}, it was also mentioned that the method developed there does not extend immediately to non-autonomous evolutionary PDEs. Here in subsection \ref{secnonauto} we deal with a non-autonomous extension of the heat equation on a compact interval with Dirichlet boundary conditions. Essentially we consider a heat equation with a continuously time varying heat diffusivity coefficient (always assumed to be positive for physical reason). This model is more physical relevant than the one studied in Section \ref{Section 2}, and for this we show that the algorithm devised in \cite{DZ} carries over with minimal modifications. 
In the final Section \ref{conclusions}, we discuss some further directions that we think would be very worthwhile to explore.

\section{The case of  a linear evolutionary PDEs of order $2N$}\label{Section 2}
In this Section, we extend the main results of \cite{DZ} to the case of a linear constant coefficient PDE of order $2N$, where in \cite{DZ} the authors restricted their analysis to the heat equation on a compact interval with Dirichlet boundary conditions.
Among the open questions posed in \cite{DZ}, it was mentioned the fact that the optimal selection of the sampling times is extremely sensitive to the  distribution of the eigenvalues of the operator that essentially controls the spatial part of the PDE. In this Section we  show that the optimal selection of sampling times {\em is essentially controlled only by order of the highest spatial derivative}, at least for the class of PDEs we consider (see Theorem \ref{Th2.10main}). So in this sense, the dependence on the spectrum is very weak in this case. To make the results directly comparable with what was obtained in \cite{DZ}, we make the strong assumption that the initial datum $f$ lives in the subspace $S\subset H^{r}_0([0,\pi])\subset L^2([0,\pi])$ consisting of functions that admit a Fourier series representation of the form $f=\sum_{k=1}^{\infty}\hat{f}_k\sin(kx)$, where $\sum_{k=1}^{\infty} k^{2r}\hat{f}^2_k<+\infty.$
Later on, we will further assume that $f\in \mathcal{F}_r$, the unit ball in $S$, defined by 
$$\mathcal{F}_r:=\{f\in S\, : \sum_{k=1}^{\infty}k^{2r}\hat{f}_k^2\leq 1\}.$$

We deal with the initial value/boundary value problem 
\beq\label{eq11}
u_t=\sum_{l=1}^N \alpha_{2l}u_{(2l)},\;0<x<\pi,\;t>0, \quad u(0,t)=u(\pi,t)=0,  \quad u(x,0)=f(x), \eeq
where $f\in S$ and $\alpha_{2l},\,l=1,2,\dots,N$ are constants. Call $\lambda: \mathbb{N}\rightarrow \mathbb{R}$, the function defined via 
\beq\label{lambdak} \lambda(k):=\sum_{l=1}^N (-1)^l \alpha_{2l}k^{2l}. \eeq
Observe that $\{\lambda(k)\}_{k\in \mathbb{N}}$ is just the spectrum of the ordinary differential operator (with respect to the $x$-variable) $L=\sum_{l=1}^N \alpha_{2l}\left(\frac{\partial}{\partial x}\right)^{2l}$ defined on $S$. 

For the problem \eqref{eq11} to be well-posed it is sufficient that for each $k\in \mathbb{N}$, $\lambda(k)\leq 0$. 
This is because the solution of the problem \eqref{eq11} is given in this case by 
\beq\label{solution1}u(x,t)=\sum_{k=1}^{\infty} \hat{f}_k e^{\lambda(k)t} \sin(kx),\eeq
where $\hat{f}_{k}$ are Fourier sine coefficients of $f(x)$. 
However, in order to extend the main results of \cite{DZ}, we will assume the following:
\beq\label{spectral1}\lambda(1)<0, \; \lambda(k+1)<\lambda(k)\; \forall k\in \mathbb{N}, \text{ and }\lim_{k\rightarrow +\infty} \lambda(k)=-\infty\eeq
In order for \eqref{spectral1} to hold, we further assume that 
\beq\label{alphaconditions} \alpha_{2l}>0 \text{ if } l \text{ is odd, and } \alpha_{2l}<0 \text{ if } l \text{ is even.} \eeq
We can immediately see how \eqref{alphaconditions} implies \eqref{spectral1}. 
\begin{proposition}\label{prop2.1}
A sufficient condition for \eqref{spectral1} to be fulfilled is that the coefficients $\alpha_{2l}$ satisfy the constraints in \eqref{alphaconditions}. 
\end{proposition}
\proof
If the conditions in \eqref{alphaconditions} are met, then $\lambda(k)$ is a polynomial in $k$ with negative coefficients, so it is clear that $\lambda(k)<0$ for each $k\in \mathbb{N}$ and $\lim_{k\rightarrow +\infty}\lambda(k)=-\infty$. Furthermore, for each $l\in \{1, \dots ,N\}$ we have $\alpha_{2l}(-1)^l (k+1)^{2l}<\alpha_{2l}(-1)^lk^{2l}$. Summing over $l$, one gets $\lambda(k+1)<\lambda(k).$
\endproof

As in \cite{DZ}, we sample at a point $x_0$ which is an algebraic number of second order, in particular we require \begin{equation}\label{algebraicnumber}|\sin(kx_0)|\geq d_0 k^{-1},\end{equation} for some $d_0>0$ and for all positive integers $k$.

\subsection{Consistency of approximation and lower bounds on optimal performance}

Our first result is the following:
\begin{theorem}\label{consitencytheorem}
Sampling $u(x,t)$ at the fixed point $x_0$ and at an increasing sequence of times $0<t_1<t_2<\dots<t_n<\dots$, either diverging to infinity or converging at some finite time, allows one to reconstruct uniquely the Fourier sine coefficients $\hat{f}_k$ and consequently $f(x)$ (in $L^2$) and $u(x,t)$. 
\end{theorem}
\proof
Introduce the following function of complex variable
$$F_0(z):=\sum_{k=1}^{\infty}c_k  z^{-\lambda(k)},$$
where the coefficients $c_k:=\hat{f}_k\sin(kx_0)$. It follows immediately that the sequence $\{c_k\}_{k\in \mathbb{N}}\in l^2$. 
Furthermore, by the current assumptions on $\lambda$, $F_0(z)$ is holomorphic in $\mathbb{D}:=\{z\in \mathbb{C}||z|<1\},$ but in general it is multi-valued there due to the fact that $\lambda(k)$ is not always an integer. We represent $z^{-\lambda(k)}$ as $\exp(-\lambda(k)\log(z))$, choose as branch cut $(-\infty, 0]$ and choose as determination of the complex logarithm $\log(z)$ the one that agrees with the real logarithm for $z\in \mathbb{R}^{+}$. In this way $F_0(z)$ becomes a single valued holomorphic function in the simply connected domain $U:=\mathbb{D}\setminus(-1,0].$
Finally observe that for each $t>0$, we have that $F_0(e^{-t})=u(x_0, t)$. Therefore, the sampling $u(x_0, t_j)=:u_j$ means receiving $F_0(z_j)$, where $z_j=e^{-t_j}.$ The sequence $\{z_j\}_{j\in \mathbb{N}}$ has a limit $z^*$. If $z^*$ is in $U$ (which means that $t_j\rightarrow t^*<\infty$), then we can invoke the identity principle:  a holomorphic function is uniquely determined by a sequence having an accumulation point on its domain of holomorphy, therefore from $F(z_j)$ we can reconstruct uniquely $c_k$, and hence $\hat{f}_k=\frac{c_k}{\sin(kx_0)}$ as $\sin(kx_0)\neq 0$ for all $k$. So the sequence $\{u_j\}_{j\in \mathbb{N}}$ uniquely determines $f\in L^2$ and consequently the solution $u(x,t)$. If on the other hand $z^*=0$ (which corresponds to the case $t_j\rightarrow +\infty$), we can not invoke the identity principle directly because the limit point is not in the interior of $U$. However we can reason as follows. Consider the a sequence of real $z_j\in U$ converging to $0\in \partial U$. Suppose $F_0(z_j)=0$ and clearly $F_0(0)=0$, since $c_0=0$. Then we want to show that $F_0$ is identically zero, namely $c_k=0$ for all $k$. This would imply the necessary identity principle in our case.   
Write $F_0(z)$ as
$$F_0(z)=z^{-\lambda(1)}\sum_{k=1}^{\infty} c_k z^{-\lambda(k)+\lambda(1)},$$
and call $g_1(z):=\sum_{k=1}^{\infty} c_k z^{-\lambda(k)+\lambda(1)}$. Since $F_0(z_j)=0$ for all $z_j$, then we have $g_1(z_j)=0$ for all $z_j$. Since $g_1$ is continuous  from the right of $0$ on the real axis, then $g_1(0)=\lim_{j\rightarrow +\infty}g_1(z_j)=0$. But this means $c_1=0$. Therefore $F_0(z)=\sum_{k=2}^{\infty} c_k z^{-\lambda(k)}$. Again  rewrite $F_0(z)=z^{-\lambda(2)}\sum_{k=2}^{\infty} c_k 
z^{-\lambda(k)}$ and call $g_2(z):=\sum_{k=2}^{\infty} c_k 
z^{-\lambda(k)}$. By the same reasoning used above one gets $g_2(0)=0$ and hence $c_2=0$. Proceeding in this way, we see that $c_k=0$ for all $k$, that is $F_0(z)=0$. This says that the identity principle can be applied also to this case, even though $0$ is not in $U$. 
Therefore, if we sample infinitely many times, the initial datum can reconstructed uniquely in $L^2$. Once the initial datum is reconstructed, then $u(x,t)$ is determined via \eqref{solution1}.
\endproof

At this point, our goal is to determine a near-optimal approximation for the initial data $f$ using only finitely many but sufficiently large number of time samplings. First we extend Theorem 3.1 from \cite{DZ} to determine optimality bounds. 
To describe bounds on optimal performance we recall few notions from the theory of manifold widths (see \cite{DVHM}). Let $\mathcal{F}_r$ be a fixed closed ball contained in $S\subset H^r_0([0\pi])\subset L^2([0,\pi])$. This ball is defined by the condition \beq\label{ballcondition}
\mathcal{F}_r:=\{f\in S :\, \sum_{k=1}^{\infty}k^{2r}|\hat{f}_k|^2\leq 1\}.
\eeq
Given $f\in \mathcal{F}_r$, we consider the problem of recovering $f$ in the $L^2$ norm and show how to modify the algorithm of \cite{DZ} to obtain a reconstruction of $f$ which is optimal in terms of rate distortion (error vs. number of measurements). We will need continuous mappings, an encoder $a: \mathcal{F}_r \rightarrow \mathbb{R}^{n}$ and a decoder $M: \mathbb{R}^{n} \rightarrow L^2$, to approximate a given $f\in \mathcal{F}_r$ as $M(a(f))$. An encoder $a$ coupled with a decoder $M$ is called a measurement algorithm. The performance of this measurement algorithm on $\mathcal{F}_r$, denoted by $\hat\delta_n(\mathcal{F}_r, L^2)$, is defined as
$$\hat\delta_n(\mathcal{F}_r, L^2)=\sup_{f\in \mathcal{F}_r}\|f-M(a(f))\|_{L^2}.$$
Since the set of all $M(y),\;y\in\mathbb{R}^n$ is an $n$-dimensional manifold, the manifold width $\delta_{n}$ is then defined as the best performance one can obtain with this scheme: 
$$\delta_{n}(\mathcal{F}_r, L^2)=\inf_{a, M}\sup_{f\in\mathcal{F}_r}\|f-M(a(f))\|_{L^2},$$
where the infimum is taken over all continuous mappings $a$ and $M$ of the above form for a fixed $n$. For $\mathcal{F}_r$ it is known that (see \cite{DVHM})
\beq\label{manifoldwidths}
\delta_n(\mathcal{F}_r, L^2)\geq c_r n^{-r}.
\eeq
We have the following extension of Theorem 3.1 of \cite{DZ}:

\begin{theorem}\label{performbound}
For any measurement algorithm (fixed or adaptive) with a continuous decoder $A_n$ we have that
$$\hat\delta_n(\mathcal{F}_r, L^2)\geq \delta_n(\mathcal{F}_r, L^2)\geq c_r n^{-r}.$$
\end{theorem}

\proof
In order to show that this bound applies also to this problem, first of all we observe that any measuring algorithm can be described by such mappings $a$ and $M$. For $a$ we take $a(f)=(x_0; t_1, \dots, t_n; u_1 \dots, u_n)$, so $a$ maps $\mathcal{F}_r$ to $\mathbb{R}^{2n+1}$ (and this covers both fixed times and adaptive times as long as the adaptive choice of times is continuous with respect to the choice of $f$). The first inequality in the statement of the Theorem is then obvious, while the second inequality is in \cite{DVHM}, for $\mathcal{F}_r$ as above, provided that the relevant map $a$ is continuous. 
 So the only thing to be checked is the fact that the encoding map $a: \mathcal{F}_r\subset L^2 \rightarrow \mathbb{R}^{2n+1}$, sending $f$ to $(x_0; t_1, \dots, t_n; u_1, \dots, u_n)$ is continuous. Since $(x_0, t_1, \dots, t_n)$ do not depend on the choice of $f\in \mathcal{F}_r$, the only claim that needs to be proved is that $u_j$ depends continuously on $f$. Let $f, g$ two functions in $\mathcal{F}_r$ and let $\hat{f}_k, \hat{g}_k$ their respective Fourier sine coefficients. 
Then \beq\label{auxiliary11}\|g-f\|_{L^2}^2=\frac{\pi}{2}\sum_{k=1}^{\infty}|\hat{g}_k-\hat{f}_k|^2.\eeq
Call $u^g$ and $u^f$ the corresponding solutions of \eqref{eq11} having as initial data $g$ and $f$ respectively. Then 
$$|u^g_j-u^f_j|=\left|\sum_{k=1}^{\infty} (\hat{g}_k-\hat{f}_k)\sin(kx_0)e^{\lambda(k)t_j} \right|\leq \sum_{k=1}^{\infty}|\hat{g}_k-\hat{f}_k|e^{\lambda(k) t_j}\leq $$
$$\leq \sum_{k=1}^{\infty}|\hat{g}_k-\hat{f}_k| e^{\lambda(k) t_1}\leq \sqrt{\frac{2}{\pi}}\|g-f\|_{L^2}\| \{e^{\lambda(k)t_1}\}_{k\in \mathbb{N}} \|_{l^2},$$
where we have used Cauchy-Schwarz inequality, \eqref{auxiliary11} and the properties of $\lambda(k)$. 
Therefore we get $|u^g_j-u^f_j|\leq C_1 \|g-f\|_{L^2},$ for a constant $C_1$ that depends only on the first time measurement $t_1$ which is always arbitrary (but greater than zero).  
\endproof

\begin{remark}
Although the authors of \cite{DZ} seem to indicate otherwise, some aspects of their method  apply   to non-autonomous linear PDEs.  For instance, consider the problem:
$$u_t=a(t)u_{xx}, \quad (t, x)\in (t_0, +\infty)\times (0, \pi), \quad u(t, 0)=u(t, \pi)=0,\quad u(t_0,x)=f(x).$$
If $a(t)$ is known, the initial time $t_0$ is known and if $\int_{t_0}^t a(s)\;ds$ is monotonic strictly increasing in $t$ (which is reasonable since if one is dealing with heat propagation then $a(t)\geq m>0$ for some $m>0$), then an adaptation of the methods developed in \cite{DZ} applies. This is because  for $t>t_0$ one has:
$$u(t, x)=\sum_{k=1}^{\infty} \hat{f}_k e^{-k^2\int_{t_0}^t a(s)\; ds} \sin(kx).$$
We will say more about this in subsection \ref{secnonauto}, where we adapt the case we analyzed in this section to this non-autonomous case. 
\end{remark}

\subsection{Time selection for near-optimal recovery}\label{optimalrecovery1}
Here we determine a sequence of times $0<t_1<t_2< \dots < t_j <\dots$ such that for sufficiently large $n$, choosing the first $n$ terms of this sequence, we can recover $f$ at the optimal rate $n^{-r}$. 

The basic idea is to use time samples $u_j:=u(x_0,t_j),\;j=1,2,\dots,n$ to create an approximation $\hat{\bar{f}}_k$ to the Fourier coefficients $\hat{f}_k$ for $k=1, \dots, n$ and then to construct the function $\bar{f}:=\sum_{k=1}^n \hat{\bar{f}}_k\sin(kx)$. The $L^2$-error of approximation to the true initial datum $f$ is given by 
\begin{equation}\label{approximationerror}
\frac{2}{\pi}\|f-\bar{f}\|^{2}_{L^2}\leq \sum_{k=1}^n |\hat{f}_k-\hat{\bar{f}}_k|^2+n^{-2r}\sum_{k\geq n+1}k^{2r}|\hat{f}_k|^2\leq  \sum_{k=1}^n |\hat{f}_k-\hat{\bar{f}}_k|^2+n^{-2r}
\end{equation}
as $f\in \mathcal{F}_r$. We need to find a sequence of times in order to approximate $\hat{f}_k$ sufficiently well so that the expression on the right side of \eqref{approximationerror} is bounded by $Cn^{-2r}$. In order to do this, let us introduce the function 
$$F(t)=\sum_{k=1}^{\infty}c_k e^{\lambda(k)t},\quad t>0$$
where $c_k:=\hat{f}_k\sin(kx_0)$. Notice that $F(t_j)=u(x_0, t_j),\;j=1,2,\dots$. First we analyze how to approximate the coefficients $c_k$ of $F(t)$ from the values $F(t_j)$, $j=1, \dots, n$. 

Let $0<t_1<t_2<\dots <t_n<\dots$ be an increasing sequence of times. Starting from $F(t_j)$, we want to derive {\em sufficient conditions} on this time sequence so that we can recover the coefficients $c_k$, $k=1,\dots, n$ with high accuracy. 

Following \cite{DZ}, we use the sample $u(x_0, t_n)$ to compute an approximation $\bar{c}_1$ of $c_1$ and then use the sample $u(x_0, t_{n-k+1})$ to compute an approximation $\bar{c}_k$ of $c_k$. 
For each $k$, we obtain $c_k$ by multiplying $F(t_{n-k+1})$ by $\exp(-\lambda(k)t_{n-k+1})$ and subtracting the remaining terms, that is, 
\begin{equation}
c_k=e^{-\lambda(k)t_{n-k+1}}F(t_{n-k+1})-\sum_{j=1}^{k-1}c_j e^{(\lambda(j)-\lambda(k))t_{n-k+1}}-\sum_{j\ge k+1}c_j e^{-(\lambda(k)-\lambda(j))t_{n-k+1}}.
\end{equation}
Now we define $\bar{c}_1:=e^{t_n}F(t_n)$ and then recursively define 
\begin{equation}\label{barc_keq}
\bar{c}_k:=e^{-\lambda(k)t_{n-k+1}}F(t_{n-k+1})-\sum_{j=1}^{k-1}\bar{c}_j e^{(\lambda(j)-\lambda(k))t_{n-k+1}}, \quad k=2, \dots, n.
\end{equation}
Then, for each $k=1,2,\dots,n$, $c_k-\bar{c}_k$ is given by
\begin{equation}\label{errorcoeff1}
c_k-\bar{c}_k=\sum_{j=1}^{k-1}(\bar{c}_j-c_j)e^{(\lambda(j)-\lambda(k))t_{n-k+1}}-\sum_{j\ge k+1} c_j e^{-(\lambda(k)-\lambda(j))t_{n-k+1}}.\end{equation}
We denote with $E_j:=|c_j-\bar{c}_j|$, the error with which we recover $c_j$ for $j\leq n$. We will first concentrate on deriving a suitable bound for $E_j$ in several subsequent lemmas, which we will use to prove the main result of this paper. 

Since $u_0\in \mathcal{F}_r$, $|c_j|^2\le |\hat{u}_j(0)|^2\le j^{-2r}\sum_{k=1}^{\infty}k^{2r}|\hat{u}_k(0)|^2\le j^{-2r}$, so we have from \eqref{errorcoeff1}
\beq\label{E1} E_1\le \sum_{j\geq 2}j^{-r}e^{-(\lambda(1)-\lambda(j))t_n}\leq 2^{-r}e^{-\delta_1 t_n}\sum_{j\geq 2}e^{-(\lambda(2)-\lambda(j))t_1}\leq A_0(t_1)e^{-\delta_1 t_n},\eeq
where $\delta_1:=\lambda(1)-\lambda(2)>0$ and $A_0(t_1)$ is a constant that depends on the initial sampling time $t_1$ and on the spectrum of the differential operator. Again use the formula \eqref{errorcoeff1} and obtain for $k\ge 2$ that
\beq\label{Ek}
E_k\leq \sum_{j=1}^{k-1}E_j e^{(\lambda(j)-\lambda(k))t_{n-k+1}}+\sum_{j\geq k+1}j^{-r}e^{-(\lambda(k)-\lambda(j))t_{n-k+1}}=:\Sigma_1(k)+\Sigma_2(k)
\eeq
We first bound $\Sigma_2(k)$. We have
\begin{equation}\label{coefferror2}
\Sigma_2(k)\leq (k+1)^{-r}e^{-\delta_k t_{n-k+1}}\sum_{j\geq k+1}e^{-(\lambda(k+1)-\lambda(j))t_{1}},
\end{equation}
where $\delta_k:=\lambda(k)-\lambda(k+1)$. Moreover, we have the following
\begin{lemma}\label{constantbound}
With the standing assumptions on the coefficients of the PDE we have
$$\sum_{j\geq k+1}e^{-(\lambda(k+1)-\lambda(j))t_1}\leq \sum_{j\geq 2}e^{-(\lambda(2)-\lambda(j))t_1}=:A_0(t_1), \quad k=2, 3,\dots.$$
\end{lemma}
\proof
The claim follows if we can show that 
$$\sum_{j=0}^{\infty}e^{-(\lambda(k+1)-\lambda(j+k+1)t_1}\leq\sum_{j=0}^{\infty}e^{-(\lambda(2)-\lambda(j+2))t_1}, \quad k=2, 3, \dots,$$
which holds provided 
$$\lambda(k+1)-\lambda(j+k+1)\geq \lambda(2)-\lambda(j+2), \quad j=0, 1, \dots,\; k=2, 3, \dots, $$
or equivalently we have for $j=0,1,\dots,\;k=2, 3, \dots$
$$\sum_{l=1}^N (-1)^l \alpha_{2l}((k+1)^{2l}-(j+k+1)^{2l})\geq \sum_{l=1}^N (-1)^l \alpha_{2l}(2^{2l}-(j+2)^{2l}).$$
Since $(-1)^l\alpha_{2l}=:\beta_l<0$ for all $l=1, \dots, N$, the last inequality holds provided 
$$\beta_{l}((k+1)^{2l}-(j+k+1)^{2l})\geq \beta_l(2^{2l}-(j+2)^{2l}), \quad l=1, \dots, N, \; k=2, 3, \dots, \; h=0,1, \dots,$$
which is true if we have
\begin{equation}\label{eqconstantbound}
(k+1)^{2l}-(j+k+1)^{2l}\leq 2^{2l}-(j+2)^{2l}, \quad l=1, \dots, N, \; k=2, 3, \dots, \; h=1, 2, \dots
\end{equation}
But then \eqref{eqconstantbound} follows from the fact that the functions $f_{k,l}(x):=(k+1)^{2l}-(x+k+1)^{2l}- 2^{2l}+(x+2)^{2l},\;l=1,2,\dots, N,\;k=2,3,\dots$ satisfy $f_{k,l}(0)=0$, and $f'_{k,l}(x)=2l\left((x+2)^{2l-1}-(x+k+1)^{2l-1} \right)<0$ for all $x\geq 0$. 
\endproof

Using Lemma \ref{constantbound}, we have that 
\begin{equation}\label{coefferror22}
\Sigma_2(k)\leq (k+1)^{-r}e^{-\delta_k t_{n-k+1}} A_0(t_1).\end{equation}

\begin{lemma}\label{inequalitypowers}
Let $a, b$ such that $1\leq a<b$. Then $b^l-a^l>b^j-a^j$ for all $l>j\geq 1$.
\end{lemma}
\proof
Clearly $1-\left( \frac{a}{b}\right)^l>1-\left( \frac{a}{b}\right)^j$ and a fortiori
$$1-\left( \frac{a}{b}\right)^l>\frac{1}{b^{l-j}}\left[1-\left( \frac{a}{b}\right)^j \right],$$
and multiplying both sides by $b^l$, we obtain the claim.
\endproof
\begin{lemma}\label{inqualitypowers2}
Let $l, m$ be positive integers with $l\geq m+1>m\geq 1$ and let $k$, $j$ be positive integers with $k-1\geq j\geq 1$. 
Then the following inequality holds:
\begin{equation}\label{eq:inequalitpowers2}
\frac{(k+1)^l-j^l}{(j+1)^l-j^l}\geq \frac{(k+1)^m-j^m}{(j+1)^m-j^m}>0.
\end{equation}
\end{lemma}
\proof 
The both ratios in \eqref{eq:inequalitpowers2} are positive as they are the ratios of positive real numbers. Therefore, to prove \eqref{eq:inequalitpowers2} under the constraints given in the statement of the lemma, it suffices to prove that 
$$\frac{(k+1)^l-j^l}{(k+1)^m-j^m}\geq \frac{(j+1)^l-j^l}{(j+1)^m-j^m}.$$
To prove this inequality, it suffices to show that the function $g_{l,m}(x):=\frac{x^l-j^l}{x^m-j^m}$ for $x\geq j+1$ and $l\geq m+1>m\geq 1$ is increasing on $x\geq j+1$. Notice that the sign of the derivative of $g_{l,m}$ is the same as the sign of the expression 
$$x^{m+l-1}(l-m)+x^{m-1}mj^l-x^{l-1}lj^m$$
which is positive if the following expression obtained from it by dividing by $x^{l-1}$ 
$$x^m(l-m)+x^m\left[m\left(\frac{j}{x}\right)^l-l\left(\frac{j}{x}\right)^m \right]$$
is positive. This is true provided $l-m>l\left(\frac{j}{x}\right)^m-m\left(\frac{j}{x}\right)^l$. To this aim, introduce the function 
$$p_{l,m}(y):=l-m+my^l-ly^m, \quad 0<y<1.$$
Now $p_{l,m}(0)=l-m>0$ and $p_{l,m}(1)=0$. Furthermore, $p'_{l,m}(j/x)=ml((j/x)^{l-1}-(j/x)^{m-1})<0$ since $l\geq m+1>m\geq 1$ and $x>j$. Therefore, the function $p_{l,m}(y)$ is positive in $(0,1)$ and thus $l-m>l\left(\frac{j}{x}\right)^m-m\left(\frac{j}{x}\right)^l$ for $x\geq j+1$. Consequently the sign of the derivative of $g_{l,m}$ is positive for $x\geq j+1$ and we are done. 
\endproof

\begin{lemma}\label{lemmarhoinequality}
The function $$g(x,y):=\frac{1}{x-y}\ln\left(\frac{(x+1)^{2N}-y^{2N}}{(y+1)^{2N}-y^{2N}}\right),$$
where $N$ is a positive integer, is positive and bounded by $2N\ln(2)$ in the domain $D:=\{(x ,y)\in \mathbb{R}^2\, : \,2\leq x, \, 1\leq y\leq x-1\}.$ 
\end{lemma}
\proof
Observe that $g$ is positive in the domain $D$. Before proving that $g$ is bounded by $\ln(4N)$ in $D$, we claim that for each fixed $y\geq 1$, $g(x,y)$ is a bounded function of $x$ and for $x\geq y+1$, it attains its maximum $h(y)$ in $D$ at $x=y+1$. To prove this claim, fix $y\geq 1$ and write $g(x,y)$ as $\frac{p_y(x)}{q_y(x)}$, where $q_y(x)=x-y$ and $p_y(x)=\ln\left(\frac{(x+1)^{2N}-y^{2N}}{(y+1)^{2N}-y^{2N}}\right).$ Then $g(x,y)$ is strictly decreasing in $x$ for $x\geq y+1$ and for each $y$ fixed provided $q_y'p_y>p_y'q_y$, that is, 
\begin{equation}\label{inequalityau}
\ln\left(\frac{(x+1)^{2N}-y^{2N}}{(y+1)^{2N}-y^{2N}}\right)>\frac{2N(x-y)(x+1)^{2N-1}}{(x+1)^{2N}-y^{2N}}.
\end{equation}
To prove this inequality, call $\lambda_1(x)$ the left hand side of \eqref{inequalityau} and $\lambda_2(x)$ the right hand side of \eqref{inequalityau}. Observe that for $x=y\geq 1$ we have $\lambda_1(y)=\lambda_2(y)$. Therefore to prove \eqref{inequalityau} it is sufficient to prove that for $x\geq y\geq 1$ one has $\lambda_1'(x)>\lambda_2'(x)$. 
Multiply both sides of $\lambda_1'(x)>\lambda_2'(x)$ by $((x+1)^{2N}-y^{2N})^2/2N$ to obtain
$$(x+1)^{2N-1}((x+1)^{2N}-y^{2N})>$$ $$\left[(x+1)^{2N-1}+(2N-1)(x-y)(x+1)^{2N-2}\right]((x+1)^{2N}-y^{2N})-2N(x-y)(x+1)^{4N-2}.$$
Adding $2N(x-y)(x+1)^{4N-2}-(x+1)^{2N-1}((x+1)^{2N}-y^{2N})$ on the both sides, the last inequality is equivalent to 
$$2N(x-y)(x+1)^{4N-2}>(2N-1)(x-y)(x+1)^{2N-2}((x+1)^{2N}-y^{2N}).$$
Dividing the both sides by $(x-y)(x+1)^{2N-2}$ and then adding $(2N-1)(x+1)^{2N}-y^{2N}$ on the both sides, the last inequality is equivalent to
$$(x+1)^{2N}+(2N-1)y^{2N}>0,$$ 
which is obviously true for $x\ge y\ge 1$. This shows that \eqref{inequalityau} holds for $x\ge y\ge 1$, from which it follows that for each fixed $y\ge 1$, $g(x,y)$ is strictly decreasing in $x$ for $x\ge y+1$ and thus has its maximum value
$$h(y)=\ln\left(\frac{(y+2)^{2N}-y^{2N}}{(y+1)^{2N}-y^{2N}}\right).$$

Next we claim that $h$ is a bounded function of $y$ for $y\geq 1$. As we notice that $h(y)$ is positive for $y\geq1$, in order to prove this claim it is enough to show that the argument of the logarithm function is bounded. Since $2^{k}\le 2^{2N}$ for all $k=1,2,\dots, 2N$, for each $y\ge 1$ we have
$$\frac{(y+2)^{2N}-y^{2N}}{(y+1)^{2N}-y^{2N}}=\frac{\sum_{k=1}^{2N}\binom{2N}{k}y^{2N-k}2^k}{\sum_{k=1}^{2N}\binom{2N}{k}y^{2N-k}}\le \frac{2^{2N}\sum_{k=1}^{2N}\binom{2N}{k}y^{2N-k}}{\sum_{k=1}^{2N}\binom{2N}{k}y^{2N-k}}=2^{2N},$$
and therefore $h(y)\leq \ln(2^{2N})=2N\ln 2$. 

Finally we have $g(x,y)\le h(y)\le 2N\ln 2$ for all pairs of $x$ and $y$ with $x\ge y+1\ge 2$. This shows that $g$ is bounded in $D$, thereby completing the proof of the lemma.  
\endproof

Now we give a choice of $t_j$ so that we can derive a bound for $\Sigma_1(k)$ comparable to the right hand side of \eqref{coefferror2}. In this way, when we combine $\Sigma_1(k)$ and $\Sigma_2(k)$ we are able to obtain the right bound for $E_k$. 

\begin{lemma}\label{samplingtimes}
Given any fixed choice of $t_1>0$, there exists $\rho>0$ such that for the choice of sequence of times $t_k:=\rho^{k-1}t_1,\; k=1,2,\dots$, we have
\begin{equation}\label{bounderrorop}
E_k\leq A_0(t_1)2^k e^{-\delta_kt_{n-k+1}}, \; k=1, 2, \dots
\end{equation}
\end{lemma}
\proof
We prove \eqref{bounderrorop} by complete induction. We see, by \eqref{E1}, that \eqref{bounderrorop} is true for $k=1$. By inductive hypothesis we have $E_j\leq A_0(t_1)2^j e^{-\delta_jt_{n-j+1}}$ for all $j<k$. Furthermore we know from \eqref{coefferror22} that $\Sigma_2(k)\leq (k+1)^{-r}A_0(t_1)e^{-\delta_k t_{n-k+1}}$, thus we just need to estimate $\Sigma_1(k)$ for this choice of sampling times. 
We thus have 
$$\Sigma_1(k)\leq \sum_{j=1}^{k-1}A_0(t_1)2^j e^{-\delta_j t_{n-j+1}}e^{(\lambda(j)-\lambda(k)) t_{n-k+1}}=A_0(t_1)\sum_{j=1}^{k-1}2^j e^{(\lambda(j)-\lambda(k)-\delta_j \rho^{k-j})t_{n-k+1}},$$
where we have used $\frac{t_{n-j+1}}{t_{n-k+1}}=\rho^{k-j}$ (See \eqref{Ek} for $\Sigma_1(k)$). Use these estimates of $\Sigma_1(k)$ and $\Sigma_2(k)$ in $\eqref{Ek}$ to obtain 
$$E_k\leq A_0e^{-\delta_k t_{n-k+1}}\left[ (k+1)^{-r}+\sum_{j=1}^{k-1}2^j e^{(\lambda(j)-\lambda(k)-\delta_j \rho^{k-j}+\delta_k)t_{n-k+1}}\right].$$
We want to show that there exists $\rho>0$ such that $\lambda(j)-\lambda(k)-\delta_j \rho^{k-j}+\delta_k\leq 0$ for $1\leq j\leq k-1$ and $k=2, 3, \dots.$
But for each such $j$ and $k$ this is equivalent to asking 
$\lambda(j)-\lambda(k+1)\leq (\lambda(j)-\lambda(j+1))\rho^{k-j}$ when we use $\delta_j=\lambda(j)-\lambda(j+1)$, and using the expression of $\lambda(.)$ from \eqref{lambdak} this is the same as proving 
$$\sum_{l=1}^{N}(-1)^l \alpha_{2l}(j^{2l}-(k+1)^{2l})\leq \left(\sum_{l=1}^N (-1)^l \alpha_{2l}(j^{2l}-(j+1)^{2l})\right)\rho^{k-j}.$$
Since $(-1)^l \alpha_{2l}<0$ for all $l=1, \dots, N$, a sufficient condition for the last inequality to hold is that 
$j^{2l}-(k+1)^{2l}\geq\left(j^{2l}-(j+1)^{2l}\right)\rho^{k-j},$ which is equivalent to 
$$\rho^{k-j}\geq \frac{(k+1)^{2l}-j^{2l}}{(j+1)^{2l}-j^{2l}}.$$ 
For all $l=1,2,\dots,N$ and for our choices of $j$ and $k$, due to Lemma \ref{inqualitypowers2} as we have $$\frac{(k+1)^{2N}-j^{2N}}{(j+1)^{2N}-j^{2N}}\geq \frac{(k+1)^{2l}-j^{2l}}{(j+1)^{2l}-j^{2l}},$$ so it is enough to find $\rho$ such that $$\rho^{k-j}\geq \frac{(k+1)^{2N}-j^{2N}}{(j+1)^{2N}-j^{2N}},$$
and taking logarithm on both sides and simplifying, this inequality is equivalent to 
$$\rho\geq \frac{1}{k-j}\ln\left(\frac{(k+1)^{2N}-j^{2N}}{(j+1)^{2N}-j^{2N}}\right)$$
but the existence of such a $\rho$ is guaranteed by Lemma \ref{lemmarhoinequality} with $\rho\ge 2N\ln 2$. 
With such a $\rho$ we have that 
$$E_k\leq A_0e^{-\delta_k t_{n-k+1}}\left[(k+1)^{-r}+\sum_{j=1}^{k-1}2^j \right]\leq A_02^k e^{-\delta_k t_{n-k+1}},$$
proving the claim. 
\endproof

Now we can prove the following Theorem extending the main result of \cite{DZ}:

\begin{theorem}\label{Th2.10main}
Consider the solution $u(x,t)$ of the IVP \eqref{eq11}. Fix $x_0$ fulfilling condition \eqref{algebraicnumber} and an arbitrary initial sampling time $t_1>0$. We then sample $u$ at $x_0$ at times $t_j:=\rho^{j-1} t_1,\;j=1,2,\dots$, where $\rho>2N\ln 2$ where $2N$ is the order of the PDE in \eqref{eq11}. Then whenever the initial data $f\in \mathcal{F}_r, \, r>0$, there exists a positive integer $n$ such that we can use the first $n$ sampled values to construct an approximation $f_n$ to $f$ that satisfies 
\begin{equation}\label{eqloweroptimalbound}
\|f-f_n\|_{L^2[0, \pi]}\leq C(r, t_1, \rho, \Delta )n^{-r}, \quad n\geq 1,
\end{equation}
where the constant $C$ can be chosen to depend only on $r$, the initial sampling time $t_1$ and the constant $\Delta:=\min_{k=1, \dots \lceil \frac{n}{2} \rceil}\{\delta(k)\}$ 
\end{theorem}

\proof
From the sampled values, we can compute the approximations $\bar{c}_k$ to $c_k:=\hat{f}_k\sin(kx_0)$ using \eqref{barc_keq}. Moreover, from Lemma \eqref{samplingtimes} we have that
$|c_k -\bar{c}_k|\leq A_0(t_1) 2^k e^{-\delta(k)t_{n-k+1}}, $ $1\leq k\leq n$. Define an approximation $\bar{f}_k:=\frac{\bar{c}_k}{\sin(kx_0)}$ to each $\hat{f}_k,\;k=1,2,\dots,n$. Then using condition \eqref{algebraicnumber}, we have that
\begin{equation}\label{errorFourier}
|\hat{f}_k-\bar{f}_k|\leq \frac{A_0(t_1) 2^k e^{-\delta(k)t_{n-k+1}}}{|\sin(kx_0)|}\leq C(t_1)k 2^k e^{-\delta(k)t_{n-k+1}},
\end{equation}
where $C(t_1):=A(t_1)/d_0$. Now we define the approximation $f_n$ to $f$ as
$f_n:=\sum_{k=1}^m \bar{f}_k \sin(kx)$, where $m:=\lceil \frac{n}{2} \rceil$. 
Then from \eqref{approximationerror}, with a constant $C_0$ depending only on $t_1$ and $r$, and defining $\Delta:=\min_{k=1, \dots m}\{\delta(k)\}$ we get
\begin{align*}
\frac{2}{\pi}\|f-f_n\|^{2}_{L^2[0, \pi]}\leq & C_0\sum_{k=1}^m k^22^{2k}e^{-2\delta(k)t_{n-k+1}}+m^{-2r}\\
\le & C_0 e^{-2\Delta t_1\rho^{\frac{n}{2}-1}}\sum_{k=1}^m e^{2\ln(k)+2k\ln 2}+m^{-2r}\\
\le & C_0 e^{-2\Delta t_1\rho^{\frac{n}{2}-1}}\sum_{k=1}^n e^{(2+2\ln 2)k}+m^{-2r}\\
\le & C_0 e^{-2\Delta t_1\rho^{\frac{n}{2}-1}} n e^{(2+2\ln 2)n}+m^{-2r},
\end{align*}
while obtaining this inequality we have used the facts that $\ln(k)\leq k$ for $k\ge 1$, $m\le n$, and $m\le \frac{n}{2}+1$ so that 
$$\min_{k=1, \dots, m}\{\delta(k)t_{n-k+1}\}=\min_{k=1, \dots, m}\{\delta(k)\rho^{n-k}t_1\}\geq \min_{k=1, \dots, m}\{\delta(k)\}\min_{k=1, \dots, m}\{\rho^{n-k}t_1\}=\Delta t_1 \rho^{\frac{n}{2}-1}.$$ 
But for a sufficiently large $n$ we have that  
$n^{2r+1}e^{(2+2\ln 2)n}e^{-2\Delta t_1\rho^{\frac{n}{2}-1}}< 1$, that is, 
$ne^{(2+2\ln 2)n}e^{-2\Delta t_1\rho^{\frac{n}{2}-1}}< n^{-2r}$
and, therefore 
$$\|f-f_n\|_{L^2[0, \pi]}\leq Cn^{-2r},$$
where $C$ is a constant depending on $t_1, r$ and $\Delta$. 
\endproof

\begin{remark} Here we remark that in general it is possible to approximate a number of Fourier coefficients greater than the number of samples. In particular, we show that with one sample we can approximate two Fourier coefficients.
Consider again the IVP
\[ u_t = u_{xx}, \;\; u(0,t) = u(\pi, t) = 0, \;\; u(x, 0) = f(x),\] with  solution $u(x, t) = \sum_{k\geq 1} \hat{f}_k e^{-k^2 t} \sin(kx)$. Suppose we have only one sample $u(x_0, t_1)$.

Let $F(t) =\sum_{k\geq 1} c_k e^{-k^2 t}$, where $c_k= \hat{f}_k  \sin(kx_0)$, and $\sin(kx_0) \neq 0$ for all $k \geq 1$. By  the procedure in \cite{DZ}, we can recover $\bar{c}_1 = e^{t_1}F(t_1)$, with error $$E_1 = | c_1 - \bar{c}_1 | \leq \frac{1}{2^r e^{3t_1} (1 - e^{-t_1})}.$$
Notice that 
\[ c_2 = e^{2^2 t_1} F(t_1) - c_1 e^{(2^2 -1^2) t_1} - \sum_{j \geq 3} c_j e^{-(j^2 - 2^2) t_1},  \]
and let $\bar{\bar{c}}_2 =  e^{2^2 t_1} F(t_1) - \bar{c}_1 e^{(2^2 -1^2) t_1}$. We use two bars instead of one to indicate that this is a different approximation compared to that found in \cite{DZ}.

Let's estimate $E_2 = |c_2 - \bar{\bar{c}}_2 |$:
\begin{align*}
E_2
\le & |c_1 - \bar{c}_1 | e^{3t_1} + \sum_{j\geq 3} |c_j| e^{-(j^2 - 2^2)t_1}\\
\le & E_1 e^{3t_1} + \sum_{j\geq 3} j^{-r} e^{-(j^2 - 2^2)t_1}\\
\le & E_1 e^{3t_1} +3^{-r} e^{-5t_1} ( 1 +e^{-7t_1} +e^{-16t_1} + \dots )\\ 
\leq & E_1 e^{3t_1} +3^{-r} e^{-5t_1} ( 1 +e^{- t_1} +e^{-2t_1} + \dots )\\
= & E_1 e^{3t_1} +3^{-r} e^{-5t_1} \frac{1}{ 1 -e^{- t_1} } \\
\leq & \frac{1}{2^r (1 - e^{-t_1})} +3^{-r} e^{-5t_1} \frac{1}{ 1 -e^{- t_1} }\\
\le &  \frac{1}{2^r (1 - e^{-t_1})}  \left ( 1 + e^{-5t_1} \right).
\end{align*}
The error $E_2$ in \cite{DZ} for $n$ samples is at most $\frac{1}{2^r (1 - e^{-t_1})} 2^2 e^{-5t_{n-1}}$, so these are comparable for $t_n$ small. 
  Notice however that in the case of the bound in \cite{DZ}, the time $t_{n-1}$ is a geometric multiple of $t_1$ which is possibly much bigger than $t_1$, so the error term $E_2$ is small if $n$ is big, while here the error term $E_2$ does not go to zero even if $t_1$ becomes very large.
 \end{remark}

 \subsection{The non-autonomous case}\label{secnonauto}
 
We address the same issues   in the case of a heat equation with a time dependent diffusivity coefficient. The more general case of a linear evolutionary PDE with time dependent coefficients is definitely much more complicated and is completely open.

Let $\alpha(t)$ be a function in $C^0([0, +\infty), \mathbb{R})$, not identically zero. For physical reasons, we assume $\alpha(t)>0$ for all $t>0$ and bounded away from zero. We study the initial/boundary value problem with an unknown $u(x,t)$ given by 
\begin{equation}\label{pde1noauto} 
u_t= \alpha(t) u_{xx},\; t>0,\;0<x<\pi, \quad u(0,t)=u(\pi,t)=0, \quad u(x,0)=f(x),
\end{equation} 
where $f\in H^r_0([0,\pi])\subset L^2$.

We set
  $\mu(k,t):= - \alpha(t)k^{2}$
and $
\lambda(k,t)=-b(t)k^2$ with $b(t):=\int_0^t \alpha(s)\, ds$. 
Since $\alpha(t)$ is positive and bounded away from zero, we have that $b(t)$ is strictly increasing, positive and not bounded from above. 
Then the solution of the problem \eqref{pde1noauto} is given by 
\beq
u(x,t)=\sum_{k=1}^{\infty}\hat{f}_k e^{\lambda(k,t)} \sin(kx)=\sum_{k=1}^{\infty}\hat{f}_k e^{-b(t)k^2} \sin(kx).
\eeq
By the very assumptions on $\alpha(t)$, the problem \eqref{pde1noauto} is well-posed. This assumptions imply $\lambda(k,t)<0$ for all $k\in \mathbb{N}$ and for all $t>0$. Moreover, again from the assumptions on $\alpha(t)$ we have that the following hold:
\beq\label{condnoauto1}
\lambda(1,t)<0\; \forall\; t>0, \quad \lambda(k+1,t)<\lambda(k,t)\; \forall\; k\in \mathbb{N}, \; \forall\; t>0,
\eeq
\beq\label{condnoauto2}
\lim_{k\rightarrow +\infty} \lambda(k,t)=-\infty \quad \text{for each fixed } t>0,
\eeq
\beq\lim_{t\rightarrow +\infty} \lambda(k,t)=-\infty \quad \text{for each fixed } k>0.\label{condnoauto3}
\eeq
 
With the choice of a new variable $T= b(t) = \int_0^t \alpha(s)\, ds$, problem \eqref{pde1noauto}  is restated as 
\[u_T=   u_{xx},\; T>0,\;0<x<\pi, \quad u(0,t)=u(\pi,t)=0, \quad u(x,0)=f(x).\]
By Lemma~\ref{samplingtimes} and Theorem~\ref{Th2.10main} (with $N=1$), the problem above  can be solved with satisfying accuracy with the choice of 
   $\rho>2\ln 2$ and for any fixed choice of $T_1>0$, as long as there  exists a (unique) increasing sequence of times $T_k,\;k=1,2,\dots$ such that $T_k \geq \rho^{k-1}T_1$. In other words, we need the choice of  $t_1>0$ such that $b(t_1)>0$,  and  an  increasing sequence $t_k,\;k=1,2,\dots$ such that $b(t_k) \geq \rho^{k-1}b(t_1)$. 
 
\begin{corollary}\label{samplingtimes1}
Consider the solution $u(x,t)$ of the IVP \eqref{pde1noauto} with the initial datum $f\in\mathcal{F}_r,\;r>0$. Fix $x_0$ such that \eqref{algebraicnumber} holds true, and  fix $\rho>2 \ln 2$. 
If we define a sequence $(t_j)_{j\geq 1} $ such that the  initial sampling time is $t_1>0$, and  $b(t_j) \geq \rho^{j-1} b(t_1)$, then 
 we have
\begin{equation}\label{bounderrorop1}
|c_k-\bar{c}_k|\leq A_0 2^ke^{-(2k+1)b(t_{n-k+1})}, \; k=1, 2, \dots, n.
\end{equation}
Thus, the constructed  approximation  $\displaystyle f_n(x):=\sum_{k=1}^{\lceil \frac{n}{2} \rceil} \bar{f}_k \sin(kx)$ ,  with $\bar{f}_k=\bar{c}_k/sin (kx_0)$, satisfies
\begin{equation}\label{eqloweroptimalbound2}
\|f-f_n\|_{L^2[0, \pi]}\leq C(r, t_1, \rho)n^{-r},
\end{equation}
where $C$ is a constant  that depends on $r$, $\rho$  and  $t_1$. 
\end{corollary}

\section{Conclusions}\label{conclusions}
We conclude mentioning some open questions that we think are definitely worthwhile exploring. 
The first one is to adapt the algorithm developed in \cite{DZ} to the case of non-autonomous linear evolutionary PDEs of types more general compared to the case we dealt with in subsection \ref{secnonauto}. For instance, a case like the one explored in Section \ref{Section 2}, but in which all the coefficients depend explicitly on time. In this set-up, even proving the unique reconstruction of the initial data with infinitely many samplings is not straightforward and it seems to require some new ideas. 

It is clear that the algorithm developed in \cite{DZ} and further investigated here is based on the fact that the PDE dynamic is equivalent (via Fourier series) to an infinite dimensional systems of ODEs that are easily integrable (in the cases analyzed here and in \cite{DZ} they are uncoupled first order linear ODEs). It would be definitely interesting to see how these ideas can be extended to {\em nonlinear integrable} PDEs, where, using a nonlinear analogue of the Fourier transform, like the Inverse Scattering Method, one can convert the PDE dynamic into the dynamic of an infinite dimensional integrable system of ODEs. One of the major problems, however, is that this integrable system of ODEs is made of ODEs that are non-trivially coupled and whose integration is not immediate, but it is based on the construction of action-angle variables.

\end{document}